\documentclass[12pt]{article}   
\usepackage{amsthm}   
\usepackage{color,graphicx}
\usepackage{amssymb}
\usepackage{bm}
\usepackage[ruled]{algorithm2e}
\usepackage[all]{xy}

\newcommand{\R}{\mathbb R}
\newcommand{\Z}{\mathbb Z}
\newcommand{\diag}{\mathrm {diag}}
\newtheorem{theorem}{Theorem}[section]
\newtheorem{proposition}{Proposition}[section]
\theoremstyle{definition}
\newtheorem{example}[theorem]{Example}
\theoremstyle{remark}

\title{An algebraic approach to Integer Portfolio problems}

\author{F. Castro \and J. Gago \and I. Hartillo \and J. Puerto \and J.M. Ucha}
\date{}
\begin{document}
\maketitle

\begin{abstract}
Integer variables allow the treatment of some portfolio optimization problems in a more realistic way and introduce the possibility of adding some natural features to the model.

We propose an algebraic approach to maximize the expected return under a given admissible level of risk measured by the covariance matrix. To reach an optimal portfolio it is an essential ingredient the computation of different test sets (via Gr\"obner basis) of linear subproblems that are used in a dual search strategy.
\end{abstract}


\section{Introduction}

Mean-variance portfolio construction lies at the heart of modern asset management and has been among the most investigated fields in the economic and financial literature. The classical Markowitz's approach, cf. \cite{Markowitz}, or \cite{Markowitz2000} for a recent reissue of his work, rests on the assumption that investors choose among $n$ risky assets to look for their corresponding weights $w_i$ in their portfolios, on the basis of
\begin{enumerate}
  \item previously estimated expected returns $\mu_i$, and
  \item the corresponding risk of the portfolio measured by the covariance matrix $\Omega$.
\end{enumerate}

Portfolios are considered {\em mean-variance efficient} in two senses:

\begin{itemize}
\item If they minimize the variance for a given admissible return $R$:
$$
\begin{array}{lll} \mbox{(MVP1)} &
\min \left( \begin{array}{ccc} w_1 & \ldots & w_n \end{array} \right) \Omega \left( \begin{array}{c} w_1  \\ \vdots \\ w_n \end{array} \right),   \\ \mbox{subject to} & \mu_1 w_1 + \cdots + \mu_n w_n \ge R, \\ & w_1 + \cdots + w_n = 1, \\ & w_j \in \R. \end{array}
$$

\item If they maximize the expected return for a given admissible risk (variance) $r^2$:
$$
\begin{array}{lll}
\mbox{(MVP2)} &
\max \mu_1 w_1 + \cdots + \mu_n w_n,  \\ \mbox{subject to} & \left( \begin{array}{ccc} w_1 & \ldots & w_n \end{array} \right) \Omega \left( \begin{array}{c} w_1 \\ \vdots \\ w_n \end{array} \right) \le r^2, \label{prob:MVP2} \\ & w_1 + \cdots + w_n = 1, \nonumber \\ & w_i \in \R.
\end{array}
$$
\end{itemize}

In problems (MVP1) and (MVP2) the weights $w_i$ stand for the percentage of a given asset in the portfolio.
It is well-known that both problems give dual views of a common analysis since they correspond to two scalarizations  of the more general bicriteria approach: obtaining the entire set $W$ of mean-variance efficient
portfolios. In this sense, (MVP1) and (MVP2) are equivalent in that $W$ can be obtained either solving (MVP1) or (MVP2) parametrically  on $R$ (admissible returns) or $r$ (admissible risks), respectively; see \cite{Sawaragi} for further details.

Although the standard statement of the mean-variance portfolio problem uses continuous variables, there are different reasons to consider integer variables, as it is pointed out in different authors \cite{Scherer:Martin, Li:Tsai, Bonami:Lejeune, Young, Corazza:Favaretto}. Our goal is to treat problem (MVP2) in its integer version based on the following considerations:

\begin{itemize}
\item In real markets, we can buy or sell only an integer amount of assets, so considering real weights in the portfolio is actually an approximation. As it is well known already for linear problems, the rounding of the real values obtained for one optimal portfolio with continuous weights may produce, in principle, an infeasible solution or a very bad approximation to the optimal integer solution. We think that this is indeed the case for a portfolio that potentially considers future contracts, as in \cite{Geman:Kharoubi} for commodity futures, to obtain lower correlations concluding on the benefits of diversification. We show this behaviour in Example \ref{ejemplo-mixta}.

\item The need to diversify the investments in a number of industrial sectors \cite{Bonami:Lejeune}.

\item The constraint of buying stocks by lots \cite{Bonami:Lejeune, Corazza:Favaretto, Jobst}.

\item Another reason for the use of integer variables, usually binary, appears in practice because portfolio managers and their clients often hate small active positions and very large number of assets, the reason being that they produce big transaction and monitoring costs. Hence it is rather usual to add the constraints associated to the following conditions:
    \begin{enumerate}
     \item there should be at least some previously decided minimum percentage (lower bounds), and
     \item there should be a maximum number of assets (upper bounds).
    \end{enumerate}
    Furthermore, transaction costs are included as linear constraints in the return equation, using decision variables, as in \cite{Li:Tsai}.

\end{itemize}

The above mentioned requirements enrich any portfolio model for real applications and require 0-1 or, more generally, integer variables. Nevertheless, as far as we know there are no specific algorithms to solve these problems. Note that methods in semi-definite programming deal with this problem, but they are oriented to the continuous case. See also \cite{Scherer:Martin} for a set of routines that handles these problems in the continuous case.
In our approach, it is natural to consider non negative integer variables $x_1, \ldots, x_n$ for the quantities of each product. Then it is necessary to take into account
\begin{itemize}
  \item the unit prices $a_1, \ldots, a_n$ of the products,
  \item the expected returns of each stock $\mu_1, \ldots, \mu_n$,
  \item and the total available budget $B$.
\end{itemize}

Then Problem (MVP2), in its integer version, can be restated as
$$\label{eq:problema-discreto}
\begin{array}{lll}
 & \max \sum_{i=1}^n \mu_i x_i = \bm{\mu}^t \bm{x} \\
 \mbox{subject to}
 & \sum_{i=1}^n a_i x_i = \bm{a}^t \bm{x} \leq B, \\
 & Q(\bm{x}) \le B^2 r^2 , \\ & \bm{x} \in (\Z_+)^n,
\end{array}
$$
where
$$
Q(\bm{x}) = \left( \begin{array}{ccc} a_1 x_1 &  \ldots &  a_n x_n \end{array} \right) \Omega \left( \begin{array}{c} a_1 x_1 \\ \vdots \\ a_n x_n \end{array} \right).
$$

\subsection{Previous works and our approach}

There have been several works with different techniques to treat Problem (MVP1) in an integer framework.

\begin{itemize}
\item For (MVP1): piecewise linear approximation in \cite{Sharpe} and \cite{Stone}, absolute deviation selection in \cite{Yamakozi:Konno}, branch-and-cut techniques in \cite{Bienstock} and minimax selection approach in \cite{Young}.

\item For (MVP1) with additional transaction cost: dual Lagrangian relaxation in \cite{Michelon:Maculan,Guignard:Kim}, linear terms transformation in \cite{Li:Chang}, separable terms transformation in \cite{Sharpe, Stone}, and parallel distributed computation reformulating the objective function into an ellipse and using piecewise linearization in \cite{Li:Tsai}.
\end{itemize}

For Problem (MVP2) the literature is not so wide, although this approach is very usual in practice for the so called {\em benchmark portfolios}, as used in \cite{Michaud:Michaud}. The results of \cite{BertsimasMS} have been in some sense a milestone in applying tools of Algebraic Geometry (namely Gr\"obner bases) to optimization, although this method is not effective for Problem (MVP2). See \cite[Ch.1,2]{Adams:Loustaunau}, \cite[Ch. 2]{Cox:Little:OShea}, or \cite{Bertsimas:Weismantel} for introductions to this subject. 

Our goal is to present a new algorithm to deal with portfolio problems with integer variables and non-linear constraints. The method is based on the computation of some test sets using Grobner bases. These bases are computed from a linear integer subproblem that contains the original linear constraints together with some new cuts induced by the non-linear constraints. The use of the \textit{reversal test-set} allows us to design
a dual search algorithm that moves from the optimal  solution of the linear subproblems towards the optimal solution of the entire portfolio problem. Our approach is new with respect to the cited references, and we will see in the last section that it is effective to deal with portfolios with number of stocks comparable to those in the literature (see \cite{Corazza:Favaretto, Li:Tsai}).


The paper is organized as follows. In Section \ref{sect:prelim} we fix the formulation of the problem to be treated with the method explained in \cite{Thomas}. In Section \ref{sect:dualsearch} the successive additions of linear constraints are explained, and the dual search algorithm based on a test-set computed from a Gr\"{o}bner basis is applied to find an optimal point of Problem (\ref{eq:problema-discreto}). It is also included an illustrative example.

Section \ref{sect:border_risk} explains the existence of a lower bound $r_b^2$ to the risk value $r^2$ below that it is not necessary to invest the whole budget to get an optimal portfolio. Section \ref{sect:computations} contains some computational experiments and Section \ref{sect:conclusions}  draws some conclusions on the paper.


\section{Preliminaries}\label{sect:prelim}
If one accepts the integer version of model (MVP2) to obtain efficient portfolios, the objective function and all constraints but one ---that is quadratic--- are linear. This is related to the form of the problem treated in \cite{Thomas}. To solve an integer programming problem $({\mathbf P})$ with linear objective function under different linear and nonlinear constraints, the following general approach can be applied (see \cite{Thomas} for a complete example):
\begin{enumerate}
\item First obtain a {\em test-set} for a linear part of $({\mathbf P})$, let us call it $({\mathbf P}')$. A test-set $T$ is a set of vectors in $\Z^n$ such that, given a feasible point $F$ of $({\mathbf P}')$, if none of the feasible points obtained adding the elements of $T$ to $F$ improves the value of the objective function, then $F$ is an optimum of $({\mathbf P}')$. A test-set of a linear integer problem can be obtained via Gr\"obner bases (\cite{Conti:Traverso}, see \cite[ch. 5]{Sturmfels} for a modern introduction).

    The best known way of obtaining these bases is using programs as {\tt 4ti2} (\cite{4ti2}), that takes advantage of the special structure of the {\em toric ideals} corresponding to linear integer problems. Programs for general Gr\"obner bases are not so good for big examples.

\item Starting at the optimum of the linear problem $({\mathbf P}')$, which is possibly an infeasible point for problem $({\mathbf P})$, use the {\em reversal} test-set ---so decreasing the objective function each time a vector of $T$ is applied--- to move throughout the set (tree) of feasible points of $({\mathbf P}')$ until you obtain feasible points {\em for the whole problem $({\mathbf P})$}. In our case, it means, portfolios with admissible risk. If this happens, one can prune the remaining feasible solutions.
\end{enumerate}

Our approach consists of applying this general idea to Problem (\ref{eq:problema-discreto}) mixing it with some bounds obtained from the continuous relaxation of the problem, to reduce the feasible region described by the linear constraints, as in \cite{Li:Tsai}.

The initial problem is
$$
\begin{array}{rl}
({\mathbf P}) & \max \, \{ \mu^t \bm{x} ~|~ \bm{a}^t \bm{x} \leq B, Q(\bm{x}) \le B^2 r^2, \bm{x} \in (\Z_+)^n \},
\end{array}
$$
and its linear relaxation is
$$
\begin{array}{rl}
({\mathbf P}') & \max \, \{ \mu^t \bm{x} ~|~ \bm{a}^t \bm{x} \leq B, \bm{x} \in (\Z_+)^n \}.
\end{array}
$$

The purpose of the following section is to explain how to obtain additional linear constraints to improve the accuracy of the linear description of problem $({\mathbf P})$, taking advantage of geometrical properties of the definition of risk. We look for some linear constraints based on the convexity of the hyperquadric given by the covariance matrix $\Omega$, which is symmetric positive definite. However, too many constraints means too many elements in the Gr\"obner basis, so the point will be to find the precise trade-off between constraints to eliminate unnecessary points in our searching, and at the same time not to increase the basis unnecessarily.

%

\section{A dual search algorithm based on a test-set}\label{sect:dualsearch}
A direct approach to problem $({\mathbf P})$ following \cite{Thomas} may lead to a non practical procedure to find the optimum. If we compute a test-set related to problem $({\mathbf P}')$, and move along the set of solutions of the linear relaxation of problem $({\mathbf P})$, the number of points to be processed is huge, even for a small number of variables. The main drawback of the procedure is the great number of integer solutions to be visited between the starting point and the feasible region of the problem $({\mathbf P})$. In order to avoid this enormous enumeration, some cuts can be added, using the convexity of the hyperquadric defined by the symmetric positive definite matrix $Q$.

We assume that a black box is available providing solutions to linear continuous optimization problems with quadratic convex constraints, as the function \texttt{fmincon} in {\sc Matlab}, the different implementations of semi-definite programming compared in \cite{Mittelmann}, or even the linear time in fixed dimension algorithm by \cite{Dyer}.

\subsection{Starting tasks}
In order to improve our representation of problem $({\mathbf P})$, we proceed as follows:
\begin{itemize}

\item The first step is the computation of the continuous solution $\bm{u}_c$ of the problem
$$
\max \, \{ \mu^t \bm{x} ~|~ \bm{a}^t \bm{x} \leq B, Q(\bm{x}) \le B^2 r^2, \bm{x} \in (\R^+)^n \},
%
$$
which gives us a return $\bm{\mu}^t \bm{u}_c = R_c$. Clearly the discrete return is less than or equal to $\lfloor R_c \rfloor$ (function \texttt{ComputeContinuousOptimum} in Algorithm \ref{DO}).

\item Secondly, we need a good discrete feasible point, which will give us a lower bound for the return. The problem $({\mathbf P})$ is always feasible, because the origin belongs to the region, but this point it is not very useful. The rounded point $\bm{u}_d = \lfloor \bm{u}_c \rfloor$ is not always feasible, as it is well-known.
 \item In order to get a feasible starting point, it is possible to decrease each coordinate until we get a feasible point. After that, the point can be improved so that the return cannot be increased in any direction inside the feasible region (function \texttt{ComputeDiscreteApprox} in Algorithm \ref{DO}).
\end{itemize}
From this point $\bm{p}_e$ we will reach the discrete optimum. Let $R_e = \bm{\mu}^t \bm{p}_e$ be the return associated with the discrete feasible point $\bm{p}_e$. A new valid formulation of the problem is
$$\label{eq:problema-2}
\max \, \{ \bm{\mu}^t \bm{x} ~|~ \bm{a}^t \bm{x} \le B, R_e \le \bm{\mu}^t \bm{x} \le \lfloor R_c \rfloor, Q(\bm{x}) \le B^2 r^2, \bm{x} \in (\Z_+)^n \}.
$$
\subsection{Addition of new linear constraints}
From the above formulation, we improve our description of problem $({\mathbf P})$ in two ways.
\begin{enumerate}
\item Adjusting hyperplanes to the hyperquadric given by upper and lower bounds on the variables.

To this end, for $j=1,\ldots,n$, we solve the continuous problems (function \texttt{ComputeLowerBounds} in Algorithm \ref{DO})
$$
\min \, \{ x_j ~|~ \bm{a}^t \bm{x} \le B, R_e \le \bm{\mu}^t \bm{x} \le \lfloor R_c \rfloor, Q(\bm{x}) \le B^2 r^2, \bm{x} \in (\R_+)^n \}.
$$
The above minimum values give us an integer lower bound $b_j$ for each variable $x_j$, applying the ceiling function. The constraints $b_j \le x_j$ are not going to be involved in the computation of the test-set through the Gr\"{o}bner basis. This is because we can write $b_{j} + y_j = x_j$, where $y_j \ge 0$, and this change of variables do not alter the coefficient matrix of the linear cuts, nor the linear cost function. Since the computation of the Gr\"{o}bner basis does depend only of this matrix, there is no extra computation time.

In a similar way, the maximization problems
$$
\max \, \{ x_j ~|~ \bm{a}^t \bm{x} \le B, R_e \le \bm{\mu}^t \bm{x} \le \lfloor R_c \rfloor, Q(\bm{x}) \le B^2 r^2, \bm{x} \in (\R_+)^n \}
$$
for $j=1,\ldots, n$, provides us upper bounds. However, these linear constraints  highly increase the size of the test-set. We will only use them to fix variables, because the upper and lower bounds of some variables are equal in many examples. This fact allows us reducing the dimensionality of the problem.

We consider the polytope
$$
P = \{ \bm{x} \in \R^n ~|~ \bm{a}^t \bm{x} \le B, \bm{\mu}^t \bm{x} \le \lfloor R_c \rfloor, \bm{\mu}^t \bm{x} \ge R_e, \bm{b} \le \bm{x} \}
$$
where $\bm{b} \le \bm{x}$ stands for the conditions $b_i \le x_i, i=1,\ldots,n$.
\item Adding nearly tangent hyperplanes to shrink the polytope.

The main idea is the addition of cuts so that the farthest regions of the polytope $P$ could be cut off. To do that, we use a point of $P$ where the function $Q$ takes its greatest value. This is equivalent to solve the continuous problem
$$\label{pmaxcont}
\max \, \{ Q(\bm{x}) ~|~ \bm{x} \in P \}.
$$
It is coded as function \texttt{ComputeMaxRisk} in Procedure \ref{NP}. Note that this problem can be efficiently solved since it is of polynomial complexity \cite{Kozlov:Tarasov:Hav}.

Let $\bm{p}_{\max}$ be a solution of Problem (\ref{pmaxcont}), $s$ be the half-line from the feasible point $\bm{p}_e$ to $\bm{p}_{\max}$, and $p' = {\mathcal Q} \cap s$ the intersection point of $s$ with the hyperquadric ${\mathcal Q}$ defined by the function $Q(\bm{x}) = B^2 r^2$.

Let $H$ be the supporting hyperplane to ${\mathcal Q}$ at the point $p'$. By the convexity of ${\mathcal Q}$, the hyperplane $H$ defines a linear half-space that contains the interior of ${\mathcal Q}$.

The coefficients of the hyperplane $H$ are real numbers, so its normal vector $\bm{n}$ may have non integer components. However, we are looking for linear constraints with integer coefficients, so we round the vector $\bm{n}$ to an integer vector $\tilde{\bm{n}} \in \Z^n$ (variable $Prec$ in Procedure \ref{NP}).
Then we proceed to compute the independent term $c$ of the tangent hyperplane to ${\mathcal Q}$ whose normal vector is equal to $\tilde{\bm{n}}$, and such that the half-space $\tilde{\bm{n}}^t \bm{x} \le \tilde{c} = \lceil c \rceil$ defines a linear half-space which contains the interior of ${\mathcal Q}$.

This process can be iterated as many times as we wish (Algorithm \ref{NP}), shrinking the polytope $P$. Nevertheless, there should be a trade-off between the number of new hyperplanes and the size of the Gr\"{o}bner basis associated with the system, so the maximum number of cuts allowed is a parameter of the algorithm. Additionally, the difference between $r_m^2 = \max \{ Q(\bm{x}) ~|~ \bm{x} \in P \}$ and the initial risk $r_0^2$ is another stopping criterion, passed as the parameter $Tol$ to the algorithm.

\end{enumerate}

\subsection{Dual iterations with the test-set}
Now after the above two phases Problem \ref{eq:problema-discreto} is transformed to
$$
\begin{array}{lll}
  & \max \bm{\mu}^t \bm{x} \label{eq:problema-3}  \\
  \mbox{ s. t. } & \bm{a}^t \bm{x} \le B, \\
  & R_e \le \bm{\mu}^t \bm{x} \le \lfloor R_c \rfloor, \\
  & \tilde{\bm{n}}_k^t \bm{x} \le \tilde{c}_k, k=1,\ldots, s,  \\
  & \bm{x} \ge \bm{b}, \bm{x}  \in (\Z_+)^n, \\
  & Q(\bm{x}) \le B^2 r^2,
\end{array}
$$
where $\tilde{\bm{n}}_k^t \bm{x} \le \tilde{c}_k, k=1,\ldots,s$ are the new cuts. The test-set $G$ is associated with the linear problem
$$
\begin{array}{lll}
  & \max \bm{\mu}^t \bm{x} \\
  \mbox{s. t.}& \bm{a}^t \bm{x} \le B, \\
  & \bm{\mu}^t \bm{x} \le \lfloor R_c \rfloor, \label{TS} \\
  & \tilde{\bm{n}}_k^t \bm{x} \le \tilde{c}_k, k=1,\ldots, s, \\
  & \bm{x} \ge \bm{b}, \bm{x} \in (\Z_+)^n.
\end{array}
$$
The condition $R_e \le \bm{\mu}^t \bm{x}$ is tested inside the tree-search, and used to prune leaves. The value of $R_e$ is updated as soon as a new feasible point with a better return is found.

Once we have the test-set of the above polytope we proceed with the resolution method.
The main bottleneck of our approach is the search over the tree defined by the test-set. The number of points to be processed is strongly related to the initial feasible point $\bm{p}_e$ found because:
\begin{enumerate}
  \item The estimated return $R_e$ defines the lowest facet of the polytope $P$ in terms of the objective value.
  \item The upper and lower bounds for the variables $x_i$ are strongly determined by $R_e$. The closer is $R_e$ to the optimal value, the narrower is the interval for each variable $x_i$.
\end{enumerate}

On the other hand, to apply the reversal test-set search we need an initial (and usually non feasible) point, but not far from feasibility. We take the point $\bm{p}_{\mathrm{bounds}}$ built by considering the independent terms of the linear constraints of Problem (\ref{eq:problema-3}), after applying the translation $b_i + y_i = x_i$, i.e.,
$$
\bm{p}_{\mathrm{bounds}} = (\bm{b}, B - \bm{a}^t \bm{b}, \lfloor R_c \rfloor - \bm{\mu}^t \bm{b}, \tilde{c}_k - \tilde{\bm{n}}^t \bm{b} )^t.
$$
 The starting point for the reversal test-set tree search is the point $\bm{p}_{\mathrm{ini}}$, the reduced of $\bm{p}_{\mathrm{bounds}}$ by the test-set $G$. This is the solution of the linear problem (\ref{TS}), as shown in \cite{Conti:Traverso}. With the reversal test-set, we search over the tree of nodes (feasible points for the linear problem) until we obtain a feasible point for the entire problem, including the quadratic constraint (function \texttt{TreeSearch} in Algorithm \ref{DO}). If the switch $SwFictBounds$ is set to true, the search is stopped in the first point that improves the estimated return given by the incumbent point $\bm{p}_e$.

Although the tree search has to end, a maximum number of processed records is passed to the procedure \texttt{TreeSearch} as a parameter. It could happen that the number of points processed exceeds the maximum allowed, and the optimum had not been found. If a new feasible point $\bm{p}'_e$ is found ($SwImprove = $ true), the bounds can be recomputed. The test-set is still valid for the new search. The only new computations are the independent terms of the hyperplanes and the reduction to find the starting point.

However, if the test-set were huge, it would be better to compute the new linear cuts given by the new estimated point $\bm{p}'_e$ and its associated Gr\"{o}bner basis. In general, the Gr\"{o}bner basis will be shorter, and the elapsed time spent in the tree search will be shortened.

\subsection{Restricted search in a region}
In the case that a new feasible point is not found after a given number of processed points ($SwImprove = $ false), it is then possible to apply a branch-and-cut technique with the bounds. Indeed, let $\bm{p}_e$ be the feasible point that gives the value $R_e$, and $\bm{b}$ the vector of lower bounds for the variables $x_i, i=1,\ldots,n$. Compute a point $\bm{b}' = \bm{b} + \alpha (\bm{p}_e - \bm{b}), 0 < \alpha < 1$ (usually $\alpha = 1/2$), which we call fictitious  bound, and consider the following problem:
\begin{equation}\label{problema-5}
\begin{array}{lll}
& \max \bm{\mu}^t \bm{x} \\
  \mbox{ s. t. } & \bm{a}^t \bm{x} \le B, \\
  & R_e \le \bm{\mu}^t \bm{x} \le \lfloor R_c \rfloor, \\
  & \tilde{\bm{n}}_k^t \bm{x} \le \tilde{c}_k, k=1,\ldots, s, \\
  & \bm{x} \ge \bm{b}', \bm{x} \in (\Z_+)^n, \\
  & Q(\bm{x}) \le B^2 r^2.
\end{array}
\end{equation}
Solving the above problem, we expect to find a new feasible point to relaunch the search process.
The idea is simple: the search is restricted to a smaller region, but the solution of the original problem is not guaranteed to be in that region. It is a heuristic technique to take advantage that this new problem does not need a new test-set. In our implementation, this search can be launched by setting the switch $SwFictBounds$ equal to true. The process is stopped as soon as a new point is found, and then we start again. If no point is found after the maximum number of allowed nodes (variable $MaxNumNodes$ in Algorithm \ref{DO}), then we stop, and the point $\bm{p}_e$ is our best value.

The pseudocode of the main algorithm is described in \texttt{DiscreteOptimum}. The procedure \texttt{NewPolytope} presents the pseudocode of the strengthening of the polytope $P$ by adding valid cuts. The switch $SwEOP$ is used to mark the end of the process.

\subsection{An illustrative example}
  Let $\bm{a} = (6075, 3105)^t$ be the vector of prices, and $\bm{\mu} = (12500, 10000)^t$ the vector of returns, with the covariance matrix equals to
  $$
  \Omega = \left( \begin{array}{rr} .832843e-4 & .485325e-4 \\ .485325e-4 & .651298e-3 \end{array} \right).
  $$
  Let $B = 9 \times 10^6$ be the budget, and $r^2 = 3 \times 10^{-5}$ the fixed risk. The continuous optimum is $\bm{u}_c = (772.754778, 215.028056)^t$, with a total return $R_c = 11809715.29$. Then $\lfloor R_c \rfloor = 11809715$, and rounding $\bm{u}_c$ we get the point $\bm{u}_d = (773, 215)$, which is not a feasible point. Subtracting from the components, we eventually reach a feasible point $\bm{p}_e = (773, 214)$, whose associated return is $R_e = 11802500$. The lower bounds $b_1$ and $b_2$ are now computed, solving first the continuous problems
  $$
  \min \, \{ x_i ~|~ \bm{a}^t \bm{x} \le B, R_e \le \bm{\mu}^t \bm{x} \le R_c, Q(\bm{x}) \le r^2 B^2, \bm{x} \in (\R^+)^2 \},
  $$
where
\begin{equation} \label{cuadrica-ejemplo1}
Q(\bm{x}) = \left( \begin{array}{cc} a_1 x_1 & a_2 x_2 \end{array} \right) \Omega \left( \begin{array}{c} a_1 x_1 \\ a_2 x_2 \end{array} \right).
\end{equation}
The respective continuous values are $\tilde{b}_1 = 752.69$, rounded to $b_1 = 753$, and $\tilde{b}_2 = 190.58$, rounded to $b_2 = 191$. We want to solve the problem
$$
\max \, \{ \bm{\mu}^t \bm{x} ~|~ \bm{a}^t \bm{x} \le B, R_e \le \bm{\mu}^t \bm{x} \le \lfloor R_c \rfloor, Q(\bm{x}) \le r^2 B^2, \bm{x} \ge \bm{b}, \bm{x} \in (\Z_+)^2 \}.
$$
In the associated linear problem
$$
\max \, \{ \bm{\mu}^t \bm{x} ~|~ \bm{\mu}^t \bm{x} + z_1 = \lfloor R_c \rfloor, \bm{a}^t \bm{x} + z_2 = B, \bm{x} \ge \bm{b}, \bm{x} \in (\Z_+)^2 \}
$$
we change the variables $x_i = y_i + b_i$, and the resulting linear problem has the same coefficient matrix. The computation of a Gr\"{o}bner basis leads to the test-set formed by vectors
$$
\begin{array}{lll}
\bm{v}_1 = (-4,5, 8775, 0)^t, & \bm{v}_2 = (-1,1, 2970, 2500)^t, & \bm{v}_3 = (0, -1, 3105, 10000)^t, \\ \bm{v}_4 = (1, -2, 135, 7500)^t, & \bm{v}_5 = (2, -3, -2835, 5000)^t, & \bm{v}_6 = (3, -4, -5805, 2500)^t.
\end{array}
$$
Now, we reduce the point
$$
\bm{p}_{\mathrm{bounds}} = (b_1, b_2, B - a_1 b_1 - a_2 b_2, \lfloor R_c \rfloor - \mu_1 b_1 - \mu_2 b_2)^t = (753, 191, 3832470, 487215)^t
$$
with the test-set to get the linear problem optimum, which is the starting point $\bm{p}_{\mathrm{ini}}$ of the tree search. In this case, $\bm{p}_{\mathrm{ini}} = (791, 192, 3598515, 2215)$. We now show the path followed by the search procedure in the tree of solutions of the linear part of the problem. In each node, we write the distance $\Delta_1$ to the continuous return associated with it, that is, $\Delta_1 = \lfloor R_c \rfloor - \bm{\mu}^t (\bm{p} + \bm{v}_i)$. The initial distance is $\Delta_e = 7215$, the difference between $\lfloor R_c \rfloor$ and $R_e$. The larger the value, the smaller the return. Therefore, values larger than $\Delta_e$ means that the corresponding branch can be pruned. The list of nodes to be processed are then ordered by $\Delta_1$. Note that black dots `$\bullet$' mean pruned nodes, and white dots `$\circ$' mean new nodes. The points are shortened to the two first components to save space.

\renewcommand\labelitemi{$\star$}
\begin{itemize}
  \item Node $\bm{p} = (791, 192)^t, \Delta_e = 7215$. Leaves $\bm{p} + \bm{v}_i, i=1,\ldots,6$:
  \begin{itemize}
    \item $\bm{p} + \bm{v}_1 = (787,197)^t, \Delta_1 = 2215, r^2 \ge r_0^2$. New node $\circ$.
    \item $\bm{p} + \bm{v}_2 = (790,193)^t, \Delta_1 = 4715, r^2 \ge r_0^2$. New node $\circ$.
    \item $\bm{p} + \bm{v}_3 = (791,191)^t, \Delta_1 = 12215 > \Delta_e$. Pruned by $\Delta_e$.
    \item $\bm{p} + \bm{v}_4 = (792, \underline{190})^t, \Delta_1 = 9715$. Pruned by bound $b_2 = 191$.
    \item $\bm{p} + \bm{v}_5 = (793, \underline{189})^t, \Delta_1 = 7215$. Pruned by bound $b_2 = 191$.
    \item $\bm{p} + \bm{v}_6 = (794, \underline{188})^t, \Delta_1 = 4715$. Pruned by bound $b_2 = 191$.
  \end{itemize}
  The above information gives rise to the following descendants.
$$
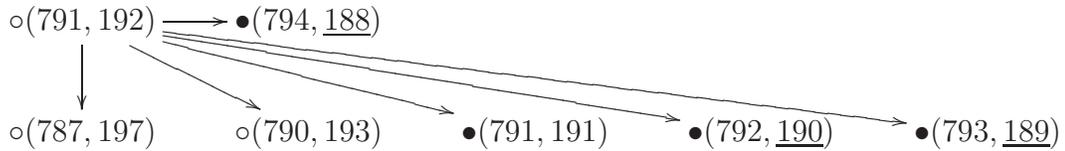
\xymatrix{
\circ (791, 192)
    \ar[d] \ar[dr] \ar[drr] \ar[drrr] \ar[drrrr] \ar[r] & \bullet (794, \underline{188}) \\ \circ (787, 197) & \circ (790,193) & \bullet (791, 191) & \bullet (792, \underline{190}) & \bullet (793, \underline{189})
}
$$
  List of ordered nodes: $\{ (787, 197)^t, (790, 193)^t \}$.

  \item Node $\bm{p} = (787, 197)^t, \Delta_e = 7215$. Leaves $\bm{p} + \bm{v}_i, i=1,\ldots,6$:
  \begin{itemize}
    \item $\bm{p} + \bm{v}_1 = (783,202)^t, \Delta_1 = 2215, r^2 \ge r_0^2$. New node $\circ$.
    \item $\bm{p} + \bm{v}_2 = (786,198)^t, \Delta_1 = 4715, r^2 \ge r_0^2$. New node $\circ$.
    \item $\bm{p} + \bm{v}_3 = (787, 196)^t, \Delta_1 = 12215 > \Delta_e$. Pruned by $\Delta_e$.
    \item $\bm{p} + \bm{v}_4 = (788, 195)^t, \Delta_1 = 9715 > \Delta_e$. Pruned by $\Delta_e$.
    \item $\bm{p} + \bm{v}_5 = (789, 194)^t, \Delta_1 = 7215 = \Delta_e$. Pruned by $\Delta_e$.
    \item $\bm{p} + \bm{v}_6 = (790, 193)^t$. Already in the list of nodes to be processed.
  \end{itemize}
  The corresponding diagram is displayed as
  $$
  
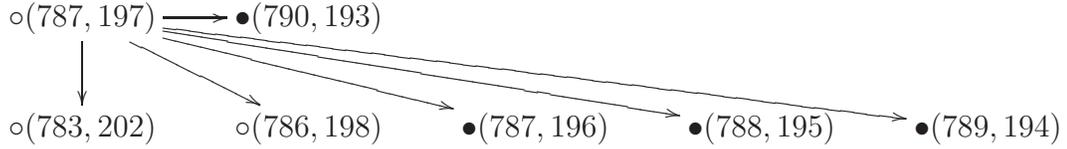
\xymatrix{
  \circ (787, 197) \ar[d] \ar[dr] \ar[drr] \ar[drrr] \ar[drrrr] \ar[r] & \bullet (790, 193) \\ \circ (783, 202) & \circ (786, 198) & \bullet (787, 196) & \bullet (788, 195) & \bullet (789, 194)
  }
  $$
  List of ordered nodes $\{ (783, 202)^t, (786, 198)^t, (790, 193)^t \}$.
  \item Node $\bm{p} = (783, 202)^t, \Delta_e = 7215$. Leaves $\bm{p} + \bm{v}_i, i=1,\ldots,6$:
  \begin{itemize}
    \item $\bm{p} + \bm{v}_1 = $\fbox{$(779, 207)^t$}, $\Delta_1 = 2215 < \Delta_e, r^2 \le r_0^2$. Feasible point, and improvement. Update $\Delta_e = 2215$.
    \item $\bm{p} + \bm{v}_2 = (782, 203)^t, \Delta_1 = 4715 > \Delta_e$. Pruned by $\Delta_e$.
    \item $\bm{p} + \bm{v}_3 = (783, 201)^t, \Delta_1 = 12215 > \Delta_e$. Pruned by $\Delta_e$.
    \item $\bm{p} + \bm{v}_4 = (784, 200)^t, \Delta_1 = 9715 > \Delta_e$. Pruned by $\Delta_e$.
    \item $\bm{p} + \bm{v}_5 = (785, 199)^t, \Delta_1 = 7215 > \Delta_e$. Pruned by $\Delta_e$.
    \item $\bm{p} + \bm{v}_6 = (786, 198)^t, \Delta_1 = 4715 > \Delta_e$. Deleted of the list of nodes to be processed.
  \end{itemize}
  The tree representation is
  $$
  
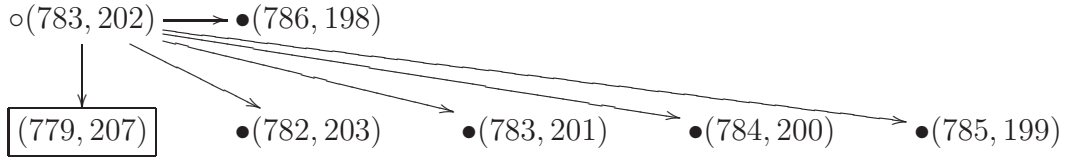
\xymatrix{
  \circ (783, 202)
 \ar[d] \ar[dr] \ar[drr] \ar[drrr] \ar[drrrr] \ar[r] & \bullet (786, 198) \\ *+[F-]{(779, 207)} & \bullet (782, 203) & \bullet (783, 201) & \bullet (784, 200) & \bullet (785, 199)
  }
  $$
  List of ordered nodes $\{ (790, 193)^t \}$.

  \item Node $\bm{p} = (790, 193)^t$. This node is pruned because $\Delta_1 > \Delta_e$, so the process is finished. The total number of processed nodes is 4.
\end{itemize}
\renewcommand\labelitemi{\textbullet}
The optimum is $(779, 207)^t$, with a difference return of $2215$ units from the continuous solution. The initial estimated point was at $7215$ units. This example shows a large difference between the rounded solution and the discrete optimum.

Although this example is very simple, we will now show the improvement that we get using the cuts. The initial polytope $P$ is defined by
$$
P = \{ \bm{x} \in \R^n ~|~  \bm{a}^t \bm{x} \le B, R_e \le \bm{\mu}^t \bm{x} \le \lfloor R_c \rfloor, \bm{x} \ge \bm{b} \}.
$$

The first maximum of the quadric $Q(\bm{x})$ described in Equation (\ref{cuadrica-ejemplo1}) on the polytope is $\bm{p}_{\max} = (753, \frac{479443}{2000})^t$ and the tangent line at the intersection point has as normal vector $(0.54452, 0.45547)^t$. Rounding to three digits, the new normal vector of the hyperplane is $(545, 455)^t$, and the independent term is $519113.7265$. Then the first cut is $545 x_1 + 455 x_2 \le 519114$. Adding this inequality to the polytope $P$, and repeating the process, we have the new maximum at point $\bm{p}_{\max} = (\frac{1979943}{2500}, 191)^t$, and the normal vector of the tangent cut is equal to $(0.56698, 0.43301)^t$. Then the new cut is $567 x_1 + 433 x_2 \le 531402$.

The linear problem, with the slack variables considered, is
$$
  \begin{array}{lll}
  & \max \bm{\mu}^t \bm{x} \\
  \mbox{s. t. } & \bm{\mu}^t \bm{x} + z_1 = \lfloor R_c \rfloor, \\
   &  \bm{a}^t \bm{x} + z_2 = B, \\
   & \bm{x} \ge \bm{b}, \\
   & 545 x_1 + 455 x_2 + z_3 = 519114, \\
   & 567 x_1 + 433 x_2 + z_4 = 531402,
  \end{array}
$$
and the test-set associated with it has $7$ elements:
$$
\begin{array}{ll}
\bm{w}_1 = (-4,5, 8775, 0, -95, 103)^t, & \bm{w}_2 = (-1,1, 2970, 2500, 90, 134)^t, \\
\bm{w}_3 = (0, -1, 3105, 10000, 455, 433)^t, & \bm{w}_4 = (1, -2, 135, 7500, 365, 299)^t, \\ \bm{w}_5 = (2, -3, -2835, 5000, 275, 165)^t, & \bm{w}_6 = (3, -4, -5805, 2500, 185, 31)^t, \\ \bm{w}_7 = ( 7, -9, -14580, 2500, 280, -72)^t.
\end{array}
$$
The starting point $\bm{p}_{\mathrm{ini}}$ is the reduction of the point
$$
\bm{p}_{\mathrm{bounds}} = ( 753, 191, 3832470, 487215, 21824, 21748)^t
$$
with the test-set. In this case, $\bm{p}_{\mathrm{ini}} = ( {\bf 779}, {\bf 207}, 3624840, 2215, 374, 78)^t$, and it is already the optimum point.

In this example, there is a large amount of money (value $3624840$) that is not invested. The explanation of this counterintuitive behaviour is because the risk $r^2$ is below a critical threshold $r_b^2$. The way to compute this threshold $r_b^2$ is the goal of the next section.

\section{A remark on admisible risks}\label{sect:border_risk}
In Problem (MVP2), there exist values of the risk $r^2$ where the optimal investment does not exhaust the available budget, i.e., the linear constraint $\bm{a}^t \bm{x} = B$ is not active in the optimal solution. Furthermore, there exists a value $r_b^2$ ({\it border risk}) below that it is not necessary to invest the whole budget to get the optimum.

The main idea is that the optimum of a linear function with a quadratic constraint $Q(\bm{x}) \le B^2 r^2$, $Q$ symmetric positive definite matrix, is found at the point on the quadric whose tangent hyperplane is parallel to the vector given by the objective function. The only problem is dealing with the negative components that the point could have.

\begin{proposition}
In Problem (MVP2), there exists a risk $r_b^2$ such that if $r^2 < r_b^2$, then the optimal investment does not need to invest the overall budget.
\end{proposition}
\begin{proof}
Given a quadric ${\mathcal C}$ with matrix
$$
{\mathcal Q} = \left( \begin{array}{cc} a_{00} & \bm{a}_0^t \\ \bm{a}_0 & Q_{00} \end{array} \right), Q_{00} \mbox{ symmetric positive definite matrix, }
$$
and $\bm{v}$ the normal vector of the hyperplane $\bm{v}^t \bm{x} = 0$, we are looking for a point $p \in {\mathcal C}$ such that the hyperplane
$$
\left( \begin{array}{cc} 1 & p^t \end{array} \right) {\mathcal Q} \left( \begin{array}{c} 1 \\ \bm{x} \end{array} \right) = 0
$$
is parallel to $\bm{v}^t \bm{x} = 0$. Then $\bm{a}_0 + Q_{00} p = \lambda \bm{v}$ for certain $\lambda \in \R$, and applying that $p$ belongs to ${\mathcal C}$ we get
$$
\lambda^2 = \frac{ \bm{a}_0^t Q_{00}^{-1} \bm{a}_0 - a_{00} }{\bm{v}^t Q_{00}^{-1} \bm{v} }, p = Q_{00}^{-1}(\lambda \bm{v} - \bm{a}_0).
$$
There are two solutions in $\lambda$, and we hold the positive one. In the case of Problem (MVP2), we have $\bm{a}_0 = \bm{0}$, $a_{00} = -r^2B^2$, and $Q_{00} = C = D \Omega D$, where $D$ is the diagonal matrix $\diag(a_1, \ldots, a_n)$. The quadric ${\mathcal C}$ is centered at the origin, the vector $\bm{v}$ is now $\bm{\mu}$, and
$$
\lambda = \frac{r B}{\sqrt{\bm{\mu}^t C^{-1} \bm{\mu} }} \mbox{ and the tangent point is } p_t =  \frac{ rB}{\bm{\mu}^t C^{-1} \bm{\mu} } C^{-1} \bm{\mu}.
$$
This point could have negative components, which means that it is outside of the feasible region of the problem. Let $J$ be the set of indexes $j$ such that the $j$-th component of vector $C^{-1} \bm{\mu}$ is positive. Let $C_J$ be the hyperquadric restricted to the intersection of hyperplanes $x_j = 0, j \in J$, and $\bm{\mu}_J$ the vector of components $\bm{\mu}_j$ with $j \in J$. The restricted hyperquadric is centered at the origin, an the optimum of the restricted problem is reached at
$$
q_t = \alpha C_J^{-1} \bm{\mu}_J, \mbox{ where } \alpha = \frac{ r B}{\sqrt{\bm{\mu}_J^t C_J^{-1} \bm{\mu}_J }}.
$$
The total amount of invested money is the dot product $q_t \bullet \bm{a}_J$, and it must be less than $B$:
$$
q_t \bullet \bm{a}_J = \frac{ r B}{ \sqrt{\bm{\mu}_J^t C_J^{-1} \bm{\mu}_J} } \bm{a}_J^t C_J^{-1} \bm{\mu}_J < B.
$$
Then
$$
r^2 < \frac{\bm{\mu}_J^t C_J^{-1} \bm{\mu}_J }{(\bm{a}_J^t C_J^{-1} \bm{\mu}_J)^2} = r_b^2.
$$
\end{proof}
It is worth noting that $r_b^2$ does not depend on the initial budget $B$.

\section{Computational results}\label{sect:computations}

This section illustrates the use of our approach in solving some integer portfolio problems with data taken from the literature. In doing that, we have implemented Algorithm \ref{DO} in {\sc Scilab}, to get a portable code, in an Intel Core2 Duo CPU, $2.53$ GHz, and 3 GB of RAM (code is available upon request for comparison purposes). The first example solves an actual problem with 44 stocks, whereas the second one shows the big sensitivity of these models with regard to the use of rounded solutions from the optimal solutions of the relaxed (continuous) formulation.
\begin{example} \label{ejemplo-eurostoxx}
This example illustrates the use of our methodology with actual data taken from  $44$ stocks indexed in Eurostoxx, from January 2003 to December 2007. The vector of initial prices is given by the prices of the stocks on January 3rd 2008 (see Table \ref{tabla-Eurostoxx}), and the returns are estimated from the monthly historical data.

\begin{table}[htbp]
\begin{tabular}{lrrlrrlrr}
Ticker & Price & Return & Ticker & Price & Return & Ticker & Price & Return \\ \hline
aca.pa & $22.7$ & $2.8$    & agn.as & $12.0$   & $0.8$    & ai.pa    & $101.6$& $12.4$ \\
alv.de & $144.9$ & $32.6$  & bas.de & $100.9$ & $24.9$  & bay.de & $61.6$ & $23.0$   \\
bbva.mc & $16.6$ &  $2.9$  & bibe.mc & $10.1$ & $2.4$   & bn.pa  & $61.2$ & $10.9$ \\
bnp.pa & $73.1$  & $11.0$    & ca.pa    & $52.2$ & $5.4$    & cs.pa  & $27.0$   & $5.6$  \\
dai.de & $64.7$  & $13.5$  & db1.de & $128.6$ & $45.8$  & dbk.de & $87.8$ & $17.8$ \\
dg.pa & $48.7$   & $14.1$  & dte.de & $14.9$ & $1.3$    & enel.mi & $8.1$ & $0.7$  \\
eni.mi & $25.1$  & $3.4$   & eoan.de & $143.8$ & $37.8$ & fora.as & $18.2$ & $1.9$ \\
fp.pa & $56.6$   & $7.8$   & fte.pa & $24.5$ & $1.5$    & g.mi    & $30.6$ & $2.2$ \\
gle.pa & $97.6$  & $15.8$  & gsz.pa & $39.5$ & $7.9$    & ing.as & $26.1$  & $5.1$ \\
isp.mi & $5.3$   & $1.3$   & lvmh.pa & $82.0$ & $13.9$    & muv2.de &   $132.0$  & $19.0$  \\
noa3.de & $25.4$ & $4.7$   & or.pa  & $96.2$ &  $9.4$   & phia.as & $28.6$ & $4.3$ \\
rep.mc & $24.9$  & $3.7$   & rno.pa & $95.2$ & $20.0$     & rwe.de  &   $95.0$   & $28.0$  \\
san.mc & $14.6$  & $3.1$   & san.pa & $62.0$   & $4.7$    & sap.de  & $34.5$ & $4.5$ \\
sgo.pa & $62.4$  & $13.0$    & sie.de & $107.1$ & $24.8$  & su.pa   & $91.1$ & $17.9$ \\
tef.mc & $21.9$  & $4.5$   & ucg.mi & $5.6$   & $0.7$ \\ \hline
\end{tabular}
\caption{Data from EuroStoxx}
\label{tabla-Eurostoxx}
\end{table}

The border risk is $r_b^2 = 0.00035$, and $B = 6000$. For the computations, all the prices and returns have been multiplied by 10 in order to work with integer variables. In this example we set the parameters to
$$
r_{\max}^2 - r_0^2 < 10^{-4}, MaxNumCuts = 4, MaxNumNodes = 10000
$$

\begin{table}[htbp]
\begin{tabular}{rrrrrrl} $r_0^2$ & $r_{\max}^2$ & cuts & basis & reduction & nodes & optimum \\ \hline
$0.0015$ &  $0.00154$ & 1 & 6657 &  & 165 &  $x_{6} = 42, x_{14} = 4, x_{16} = 29,$ \\ & & & (9 s.) & (22 s.) & (99 s.) & $x_{20} = 5, x_{28} = 1, x_{36} = 8$. \\ \hline
$0.0020$ &  $0.00205$ & 1 & 40256 &  & 137 & $x_6 = 51, x_{14} = 5, x_{16} = 19, $ \\ & & & (446 s.) & (240 s.) & (497 s.) & $x_{20} = 1, x_{28} = 1, x_{36} = 12$. \\ \hline
$0.0025$ &  $0.00256$ & 2 & 27082 &  & 32  & $x_{6} = 59, x_{14} = 6, x_{16} = 13,$ \\ & & & (194 s.) & (421 s.) & (83 s.) & $x_{28} = 2, x_{36} = 10$. \\ \hline
$0.0030$ &  $0.00301$ & 1 & 12504 &  & 62 & $x_{6} = 64, x_8 = 1, x_{14} = 8, $ \\ & & & (26 s.) & (142 s.) & (81 s.) & $x_{16} = 1, x_{28} = 1, x_{36} = 10, x_{37} = 1$ \\ \hline
$0.0035$ & $0.00351$ & 1  & 2357  &  & 0  & $x_6 = 68, x_{14} = 10, x_{16} = 1, x_{36} = 5$. \\ & & & (1 s.) & (4 s.) & (0 s.) \\ \hline
$0.0040$ & $0.00430$ & 0  & 569   &  & 9844 & $x_{6} = 74, x_{14} = 10, x_{16} = 1,$ \\ & & & (0 s.) & (6 s.) & (821 s.) & $x_{28} = 2, x_{36} = 1$. \\
         & $0.00404$ & 1  & 11924 &  & 11 & $x_{6} = 74, x_{14} = 10, x_{16} = 1,$ \\ & & & (32 s.) & (121 s.) & (10 s.) & $x_{28} = 2, x_{36} = 1$. \\ \hline
$0.0045$ & $0.00451$ & 1  & 7087 &  & 0   & $x_{6} = 84, x_{14} = 6, x_{16} = 1,$ \\ & & & (14 s.) & (33 s.) & (0 s.) & $x_{28} = 1$. \\ \hline
$0.0050$ & $0.00508$ & 0  & 357  &  & 6   & $x_{6} = 91, x_{14} = 3, x_{28} = 1$. \\ & & & (1 s.) & (0 s.) & (0 s.) \\ \hline
\end{tabular}
\caption{Discrete optimums for different values of risk $r_0^2$}
\label{DiscreteOpt-01}
\end{table}

In Table \ref{DiscreteOpt-01} we show the results of the application of our algorithm assuming different risk levels. Each risk level is organized in a block of two rows. The first one gives the corresponding information, and below we report on the elapsed time to obtain these elements. Column $r^2_{\max}$ contains the greatest risk reached in the polytope with the number of added tangents as shown by the following column. Column `basis' denotes the number of elements of the computed test-set, and column `processed' contains the number of new nodes found and explored. Finally, column `optimum' is a feasible point with the best return. The time in seconds of these tasks appears in parenthesis under the columns `basis', `reduction' and `nodes', respectively.

It is worth remarking the case $r^2 = 0.0040$. The first row shows the number of processed nodes until an optimal point was reached, with no added cutting hyperplane. The second row gives us an example of the effectiveness of adding new cuts. The test-set is computed very quickly, although the number of elements is big. However, the number of processed points is very small, and hence it is also small the total elapsed time.

\begin{table}[htbp]
\begin{tabular}{rrrrrrrl} $r^2_{\mathrm{max}} $ & tang. & basis & reduction & $Sw 1$ & $Sw 2$ & nodes & improvement \\ \hline
$0.00118$ & 3 & 32495 & & 0 & 1 & max & \\ & & (244 s.) & (333 s.) & & & (426 s.) \\ 
 &  &  & & 1 & 0 & 88 & $x_6 = 35, x_8 = 36, x_{14} = 2, $ \\ & & & (176 s.) & & & (0 s.) & $x_{16} = 27, x_{20} = 8, x_{28} = 31,$ \\   & & & & & & & $x_{35} = 3, x_{36} = 3, x_{43} = 1$ \\
 \hline
$0.00107$ & 4 & 16930 & & 0 & 1 & max & $x_6 = 34, x_8 = 22, x_{14} = 2, $ \\ & & (52 s.) & (130 s.) & & & (393 s.) & $x_{16} = 29, x_{20} = 9, x_{28} = 24,$ \\  & & & & & & & $x_{35} = 3, x_{36} = 3, x_{43} = 1 $ \\ & & & & & & & reached at 2894 nodes in 48 s. \\
\hline
$0.00107$ & 4 & 18637 & &  0 & 1 & max & \\ & & (49 s.) & (114 s.) & & & (439 s.) \\ 
 & & & & 1 & 0 & 9759 & $x_6 = 33, x_8 = 29, x_{14} = 2, $ \\ & & & (80 s.) & & & (785 s.) & $x_{16} = 31, x_{20} = 9, x_{28} = 26,$ \\ & & & & & & & $x_{35} = 2, x_{36} = 3$ \\ \hline
$0.00105$ & 4 & 14670 & & 0 & 1 & max & $x_6 = 34, x_8 = 32, x_{14} = 2, $ \\ & & (28 s.) & (79 s.) & & & (627 s.) & $x_{16} = 29, x_{20} = 8, x_{28} = 10,$ \\  & & & & & & & $x_{35} = 3, x_{36} = 4, x_{43} = 2$ \\ & & & & & & & reached at 1680 nodes in 29 s. \\ \hline
$0.00101$ & 2 & 2613 & & 0 & 0 & 2648 & $x_6 = 34, x_8 = 32, x_{14} = 2, $ \\ & & (1 s.) & (10 s.) & & & (853 s.) & $x_{16} = 29, x_{20} = 8, x_{28} = 10,$ \\   & & & & & & & $x_{35} = 3, x_{36} = 4, x_{43} = 2$ \\ & & & & & & & {\bf optimum} \\ \hline
\end{tabular}
\caption{Steps in the computation for $r_0^2 = 0.0010$}
\label{DiscreteOpt-02}
\end{table}

 Table \ref{DiscreteOpt-02} contains all the iterations done in the computation of the optimum for the case $r_0^2 = 0.0010$. The column $Sw 1$ refers to the variable $SwFictBounds$, and $Sw 2$ to $SwNumNodes$. The first one is true when fictitious  bounds are used, and the second one is true when the number of processed records is greater than $MaxNumNodes = 10000$. The column `improvement' contains a new point found with better return than the initial point. Again, the time in seconds of each task appears in parenthesis.

\end{example}

\begin{example} \label{ejemplo-mixta}
This example is devoted to show the difference between the rounded continuous solution and the integer solution of a portfolio problem. It consists of a mixing of usual stocks (Microsoft and General Electric) with the value of a future contract based on oil, as in \cite{Chng}, so the number of values is $n = 3$. The initial data is given by

\begin{center}

\begin{tabular}{lrr}
 Stock & Price $(a_i)$ & Return $(\mu_i)$ \\ \hline
 MSFT  & $35.22$ & $3.64$ \\
 GE    & $36.76$ & $3.64$ \\
 Oil   & $4000$  & $10000$ \\ \hline
\end{tabular}

\end{center}

and the covariance matrix is
$$
\Omega =
\left( \begin{array}{rrr} 0.003250634   & 0.000654331 & 0.022513263 \\
0.000654331 & 0.001578359 & -0.006610861 \\
0.022513263 & -0.006610861 &    26.35846804 \end{array} \right).
$$

We fix the risk to $r_0^2 = 1.52$, and compute the optimum for different budgets $B$. The test-set associated with the constraints
$$
\bm{a}^t \bm{x} \le B,  R_e \le \bm{\mu}^t \bm{x} \le \lfloor R_c \rfloor, \bm{x} \in (\Z_+)^n
$$
has $2663$ elements. The basis remains equal for all the considered cases. However, the capacity of computation is run out when only one more cut is added. If we take as initial point $\bm{p}_e$, the discrete approximation given by rounding, the tree searching was unable to reach the optimal point for $MaxNumNodes = 50000$. This fact is reported as `E' in Table \ref{tabla-mixta}

The function \texttt{ComputeDiscreteApprox} should be changed to get a better discrete approximation than the rounded value. If we take an approximation based on the best value for the future contract, we get Table \ref{tabla-mixta}.

\begin{table}[htbp]
\begin{tabular}{rrrrrr} & continuous & \\ Budget & optimum & return & nodes & optimum & $R_d$ \\ \hline
50000 & $(1079.87, 0, 2.99)$ & $33848.34$ \\
 & discrete approx. & \\
 & $(1192, 0, 2)$ & \fbox{$24338.88$}&  E & \\
 & $(219, 824, 3)$ & $33796.52$& 6066 & $(314, 705, 3)$ & $\mathbf{33815.06}$ \\ \hline

 75000 & $(1619.80, 0, 4.49)$ & $50772.52$ \\
 & discrete approx. \\
 & $(1675, 0, 4)$ & $46097.00$ & 22790 & $(1675, 0, 4)$ & $46097.00$ \\ \hline

 100000 & $(2159.74, 0, 5.98)$ & $67696.69$ \\
 & discrete approx. \\
 & $(2271, 0,5)$ & \fbox{$58266.44$} & E & \\
 & $(439, 1646, 6)$ & $67596.69$ & 22991 & $(687, 1409, 6)$ & $\mathbf{67629.44}$ \\ \hline
\end{tabular}
\caption{Mixed example}
\label{tabla-mixta}
\end{table}

It is easy to see the enormous difference between the return of the discrete approximation and the corresponding return of the discrete optimum for each budget.

\end{example}
\section{Conclusions}\label{sect:conclusions}
We have presented an algorithm to deal with portfolio problems with integer variables and non-linear constraints. The presented model was not previously treated as far as we know. The method is based on the computation of some test sets using Gr\"{o}bner bases, an algebraic tool. These Gr\"{o}bner bases are computed from a linear integer subproblem that contains the original linear constraints and some new cuts induced by the non linear constraints. The reversal test-set, given by the Gr\"{o}bner basis, allows us to perform a dual search algorithm from the optimal solution of the linear subproblem towards the optimal solution of the whole portfolio problem. This technique has allowed us to solve problems of size similar to the exposed in \cite{Corazza:Favaretto, Li:Tsai}.

\section*{Acknowledgments}
This research was supported by grants P06-FQM-01366  and FQM-333 (Junta de Andaluc\'{\i}a), MTM2007-64509 and MTM2007-67433-C02-01 (Ministerio de Educaci\'{o}n y Ciencia). We wish to thank Javier Nogales for beneficial discussions.

\vspace{1cm}

\noindent F. Castro Jim\'{e}nez. Dpto. de \'{A}lgebra, Univ. de Sevilla,
Apdo. 1160, 41080 Sevilla, e-mail:\texttt{ castro@us.es} \\

\noindent J. Gago Vargas. Dpto. de \'{A}lgebra, Univ. de Sevilla,
Apdo. 1160, 41080 Sevilla, e-mail:\texttt{ gago@us.es} \\

\noindent I. Hartillo Hermoso. Dpto. de Matem\'{a}tica Aplicada I, Univ. de Sevilla, E.T.S. de Ingenier\'{\i}a Inform\'{a}tica, Avda. de Reina Mercedes, s/n, 41012 Sevilla, Spain, e-mail:\texttt{ hartillo@us.es} \\

\noindent J. Puerto Albandoz, Dpto. de Estad\'{\i}stica e Investigaci\'{o}n Operativa, Univ. de Sevilla, Facultad de Matem\'{a}ticas, C/ Tarfia, s/n, 41012 Sevilla, e-mail:\texttt{ puerto@us.es} \\

\noindent J.M. Ucha Enr\'{\i}quez. Dpto. de \'{A}lgebra, Univ. de Sevilla,
Apdo. 1160, 41080 Sevilla, e-mail:\texttt{ ucha@us.es}

\newpage
\section*{Appendix}
\begin{procedure}
  \SetKwFunction{TangentToQuadric}{TangentToQuadric}
  \SetKwFunction{TangentToQuadricV}{TangentToQuadricV}
  \SetKwFunction{ComputeMaxRisk}{ComputeMaxRisk}
  \SetKwFunction{Round}{Round}
  \SetKwFunction{Ceil}{Ceil}
  \SetKwFunction{Polytope}{Polytope}
  $NumCuts = 0$ \;
  $\bm{p}_{\mathrm{max}}, r_m^2 = \ComputeMaxRisk(P, Q)$ \;
  \tcc*[f]{solve the problem $\max Q(\bm{x}), \bm{x} \in P$.}\;
  \While{$r_m^2 - r_0^2 > Tol$ {\bf and} $NumCuts \le MaxNumCuts$}
   {
        $s = {\bm p}_e + \lambda (\bm{p}_{\mathrm{max}} - {\bm p}_e), \lambda \ge 0$ \;
        ${\bm p}' = s \cap Q$ \;
        $\bm{v} = \TangentToQuadric(Q, \bm{p}')$ \;
        $DirApprox = \Round(\bm{v}, Prec)$ \; \tcc*[f]{round with number of digits = $Prec$}\;
        $Coef = \TangentToQuadricV(DirApprox, Q)$ \;
        \tcc*[f]{independent term of the tangent hyperplane $Q$ and normal vector equal to $DirApprox$} \;
        $Coef = \Ceil(Coef)$ \;
        \tcc*[f]{the best integer to leave the quadric in a half-space} \;
        $H:= DirApprox^t \bm{x} - Coef$ \;
        \tcc*[f]{new linear cut} \;
        $NumCuts = NumCuts  + 1 $\;
        $P = \Polytope(P, H)$ \; \tcc*[f]{ add a new cut to polytope $P$ } \;
        $\bm{p}_{\mathrm{max}}, r_m^2 = \ComputeMaxRisk(P, Q)$ \;
   }
   \Return{P}
\caption{NewPolytope(polytope $P$, matrix $Q$,  point $\bm{p}_e$, $Tol$, $r_0$, $MaxNumCuts$) \label{NP}}
\end{procedure}

\begin{algorithm}
  \SetKwFunction{ComputeContinuousOptimum}{ComputeContinuousOptimum}
  \SetKwFunction{ComputeDiscreteApprox}{ComputeDiscreteApprox}
  \SetKwFunction{ComputeLowerBounds}{ComputeLowerBounds}
  \SetKwFunction{Polytope}{Polytope}
  \SetKwFunction{ComputeMaxRisk}{ComputeMaxRisk}
  \SetKwFunction{TangentToConic}{TangentToConic}
  \SetKwFunction{NewPolytope}{NewPolytope}
  \SetKwFunction{ComputeMaxRisk}{ComputeMaxRisk}
  \SetKwFunction{ComputeTestSet}{ComputeTestSet}
  \SetKwFunction{Reduce}{Reduce}
  \SetKwFunction{TreeSearch}{TreeSearch}
  \KwData{budget $B$, risk $r_0^2$, matrix $Q$, vector $\bm{a}$, vector $\bm{\mu}$, $MaxNumCuts$, $MaxNumNodes, Tol$}
  \KwResult{$Optimum$, $NumNodesProc$}
    $\bm{p}_c = \ComputeContinuousOptimum(B, r_0^2, Q, \bm{a}, \bm{\mu}), g_c = \bm{\mu}^t \bm{p}_c, \alpha = 1/2$ \;
    $\bm{p}_e, g_e = \ComputeDiscreteApprox(\bm{p}_c, B, \bm{a}, r_0^2, Q)$ \;
    $SwEOP = \mathrm{false}, SwFictBounds = \mathrm{false}, SwImprove = \mathrm{false}, SwNumNodes = \mathrm{false}, ListOfVariables = (1: \dim)$ \;
    \While{{\bf not} $SwEOP$}
    {
        \If{{\bf not} $SwFictBounds$}
            {
            $\bm{b}, ListOfVariables = \ComputeLowerBounds(ListOfVariables, \bm{p}_e)$ \;
            }
        $P = \Polytope(\bm{a}^t \bm{x} \le B, \bm{\mu}^t \bm{x} \le g_c, \bm{\mu}^t \bm{x} \ge g_d, \bm{x} \ge \bm{b} )$ \;
        $NumNodesProc = 0$ \;
        \While{{\bf not} ( $SwEOP$ {\bf or} $SwImprove$ {\bf or} $SwNumNodes$)}
        {
            $P = \NewPolytope(P, Q, \bm{p}_e, Tol, r_0, MaxNumCuts)$ \;
            $G = \ComputeTestSet(P), \bm{p}_{\mathrm{ini}} = \Reduce(\bm{p}_{\mathrm{bounds}}, G)$ \;
            $SwNumNodes, SwImprove, Optimum = \TreeSearch(\bm{p}_{\mathrm{ini}}, G, Q, \bm{a}, B, r_0^2, \bm{b}, \bm{p}_e, MaxNumNodes, SwFictBounds)$ \;
            \eIf{$SwFictBounds$}{
                \eIf{{\bf not} $SwImprove$}{
                    $SwEOP = \mathrm{true}$ \;
                    }
                    {
                    $\bm{p}_e = Optimum, g_e = \bm{\mu}^t Optimum, \, SwFictBounds = \mathrm{false}$ \;
                    }
            }
            {
            \eIf{{\bf not} $SwNumNodes$}
            {
                    $SwEOP = \mathrm{true}$ \;
            }
            {
                \If{ {\bf not} $SwImprove$}
                    {
                    $\bm{b} = \bm{b} + \alpha (\bm{p}_e - \bm{b}), \, SwFictBounds = \mathrm{true}$ \;
                    }

            }
            }
        }
    }
\caption{DiscreteOptimum}
\label{DO}
\end{algorithm}

\end{document}